\documentclass[12pt]{article}

\usepackage{color}
\usepackage{graphicx}

\def\defeq{:=}
\def\eqdef{=:}

\newcommand{\mylabel}[1]{\label{#1}}
\newcommand{\myref}[1]{\ref{#1}}
\newcommand{\myeqref}[1]{\eqref{#1}}
\newcommand{\myeqlab}[1]{\label{#1}}
\newcommand{\eqref}[1]{(\ref{#1})}
\newcommand{\operatorname}[1]{{\rm#1}~}
\newcommand{\text}[1]{{\rm#1}}

\renewcommand{\Re}{\operatorname{Re}}
\renewcommand{\Im}{\operatorname{Im}}

\newcommand{\Ck}[1]{\spC^{#1}}

\newcommand{\Cone}{\Ck1}

  \newcommand{\acf}{a}

\newcommand{\mathds}[1]{\mathbf{#1}}

\newcommand{\Z}{\mathds{Z}}
\newcommand{\defm}[1]{\emph{#1}}
\newcommand{\subeq}[2]{\mathord{\underbrace{\mathop{#1}}_{#2}}}

\newcommand{\graph}{\operatorname{graph}}
\newcommand{\sign}{\operatorname{sign}}

\newcommand{\spC}{\mathcal{C}}

\def\XXint#1#2#3{{\setbox0=\hbox{$#1{#2#3}{\int}$}
\vcenter{\hbox{$#2#3$}}\kern-.5\wd0}}

\def\XXsum#1#2#3{{\setbox0=\hbox{$#1{#2#3}{\sum}$}
\vcenter{\hbox{$#2#3$}}\kern-.5\wd0}}

\newcommand{\vv}{\vec v}

\newcommand{\aint}{{-\kern-5.3mm\int}}

\newcommand{\half}{\frac12}

\newcommand{\cli}[2]{{[#1,#2]}}

\newcommand{\vlen}[1]{|#1|}

\newcommand{\csep}{\quad,\quad}

\newcommand{\pt}{\partial_t}

\newcommand{\dotp}{\cdot}
\newcommand{\crossp}{\times}

\newcommand{\ndiv}{\nabla\dotp}
\newcommand{\ncurl}{\nabla\crossp}
\newcommand{\nperp}{\nabla^\perp}

\newcommand{\isoba}[1]{\mathcal{H}_0^{#1}}

\newcommand{\conv}{\rightarrow}

\newcommand{\Mor}[2]{\operatorname{Mor}}

\newcommand{\const}{\text{const}}

\newcommand{\set}[1]{\{#1\}}

\newcommand{\Lap}{\Delta}

\newcounter{cons}

\newcommand{\pd}[1]{\partial_{#1}}

\newcommand{\rhsf}{F}

\newcommand{\rmin}{\underline\rad}
\newcommand{\Pola}{\Theta}

\newcommand{\Ber}{B}
\newcommand{\isenc}{\gamma}

\newcommand{\polalo}{\pola_0}
\newcommand{\polahi}{\pola_1}

\newcommand{\www}{w}
\newcommand{\wwi}{w_\infty}
\newcommand{\wpot}{\Phi}

\newcommand{\Dom}{\Omega}

\newcommand{\Body}{B}
\newcommand{\vx}{v^x}
\newcommand{\vy}{v^y}

\newcommand{\cz}{\overline z}
\newcommand{\pcz}{\pd\cz}
\newcommand{\pz}{\pd z}

\newcommand{\Gam}{\Gamma}

\newcommand{\vcav}{v_*}

\newcommand{\vort}{\omega}
\renewcommand{\Re}{\operatorname{Re}}
\renewcommand{\Im}{\operatorname{Im}}

\newcommand{\dens}{\varrho}
\newcommand{\vpot}{\phi}

\newcommand{\stf}{\psi}

\newcommand{\piv}{p}
\newcommand{\pif}{\hat\piv}
\newcommand{\ipif}{\pif^{-1}}
\newcommand{\pp}{P}
\newcommand{\ppf}{\hat\pp}

\renewcommand{\vec}[1]{\mathbf{#1}}

\newcommand{\rad}{r}
\newcommand{\vn}{\vec n}
\newcommand{\vs}{\vec s}
\newcommand{\xx}{{\vec x}}

\newcommand{\pola}{\theta}
\newcommand{\ssnd}{c}
\newcommand{\Mach}{M}

\newcommand{\vvi}{\vv_\infty}

\newcommand{\hdiv}{\hat\tau}

\newcommand{\vms}{\mu}
\newcommand{\vmsmax}{\overline\vms}

\title{Nonexistence of irrotational flow around solids with protruding corners}
\author{Volker Elling}
\date{}

\begin{document}

\maketitle

\abstract{%
  We motivate and discuss several recent results on non-existence of irrotational inviscid flow around bounded solids 
  that have two or more protruding corners, complementing classical results for the case of a single protruding corner.
  For a class of two-corner bodies including non-horizontal flat plates, compressible subsonic flows do not exist.
  Regarding three or more corners, bounded simple polygons do not admit compressible flows with arbitrarily small 
  Mach number, and any incompressible flow has unbounded velocity at at least one corner.
  Finally, irrotational flow around smooth protruding corners with non-vanishing velocity at infinity does not exist. 
  This can be considered vorticity generating by a slip-condition solid in absence of viscosity. 
}


\begin{figure}
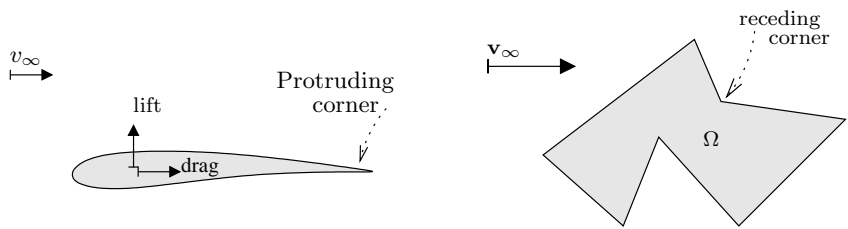
%
  \parbox{2.5in}{\input{kj.pstex_t}}\hfill\parbox{3in}{\input{polygonflow.pstex_t}}%
  \caption{Left: flow around a body that is smooth except for one protruding corner;
    right: flow around a polygon.}%
  \label{fig:onecorner}%
\end{figure}%

\section{Equations}

Consider inviscid flow in the 2d plane around a solid body whose boundary is smooth except for one or more corners. 
This setup is an old and important problem in fluid dynamics. 
The case of a single corner has received particular attention, since it idealizes cross-sections of aircraft wings; 
the corresponding Kutta-Joukowsky theory is the basis for the understanding of lift, the upward force on a horizontally moving body. 

To explain the problem of corners, we first recall some properties of compressible flow: 
\begin{eqnarray} 
  0 &=& \pt\dens + \ndiv(\dens\vv) \quad,\quad
  0 = \pt\vv + \vv\dotp\nabla\vv + \nabla\piv 
  \myeqlab{eq:unsteady}
\end{eqnarray} 
where $\vv$ is velocity, $\dens$ density; $\piv=\pif(\dens)$ is defined up to an additive constant by
\begin{eqnarray} \pif_\dens &=& \dens^{-1}\ppf_\dens \end{eqnarray} 
where $\pp=\ppf(\dens)$ is \defm{pressure}. We consider the \defm{polytropic} pressure law
\begin{eqnarray} \pp &=& \ppf(\dens) = \dens^\isenc \end{eqnarray} 
with \defm{isentropic coefficient} $\isenc>1$, so that
\begin{eqnarray} \pif(\dens) &=& \frac{\isenc}{\isenc-1} \dens^{\isenc-1}. \myeqlab{eq:pif}\end{eqnarray} 

Linearizing \myeqref{eq:unsteady} around a $\vv=0$ and $\dens=\overline\dens=\const>0$ background yields
\begin{eqnarray} 0 &=& \pt\dens + \overline\dens~\ndiv\vv \csep 0 = \pt\vv + \overline\dens^{-1}\ppf_\dens(\overline\dens)\nabla\dens \end{eqnarray} 
$\pt$ of the first minus $\overline\dens~\ndiv$ of the second equation yields
\begin{eqnarray} 0 &=& \pt^2\dens - \ppf_\dens(\overline\dens)\Lap\dens \end{eqnarray} 
which is the linear wave equation; this motivates defining the \defm{speed of sound}
\begin{eqnarray} \ssnd &=& \sqrt{\ppf_\dens(\dens)}. \end{eqnarray} 

For steady flow
\begin{eqnarray} 
  0 &=& \ndiv(\dens\vv), \myeqlab{eq:masssteady} \\
  0 &=& \vv\dotp\nabla\vv + \nabla\piv,  \myeqlab{eq:vsteady}
\end{eqnarray} 
the dot product of \myeqref{eq:vsteady} with $\vv$ yields
\begin{eqnarray} 0 &=& \vv\dotp\nabla\vv\dotp\vv + \vv\dotp\nabla\piv = \vv\dotp\nabla(\half\vlen\vv^2+\piv) \end{eqnarray} 
Hence the \defm{Bernoulli constant}
\begin{eqnarray} \Ber &=& \half\vlen\vv^2 + \pif(\dens) \end{eqnarray} 
is constant along \defm{streamlines}, i.e.\ integral curves of $\vv$; this is the \defm{Bernoulli relation}. 
On each streamline we can solve for
\begin{eqnarray} \dens = \ipif(\Ber-\half\vlen\vv^2). \myeqlab{eq:dens-pifi} \end{eqnarray} 
$\ipif$ is generally undefined for negative arguments
due to fractional exponents (see \myeqref{eq:pif}), 
in particular when $\isenc=\frac53$ (helium and other noble gases) or $\isenc=\frac75$ (air).
Hence $|\vv|$ may not exceed the \defm{limit speed} $\vcav=\sqrt{2\Ber}$; 
at the limit speed the density reaches $0$. There is no meaningful way to extend the model to higher speeds,
as becomes clear by putting the observation into a more physical form:
groups of gas particles cannot acquire arbitrarily high speed by moving to regions of increasingly lower pressure.
(A related observation: gas inside a piston expanding to near-vacuum cannot perform unbounded mechanical work, 
or put differently, it takes only a finite amount of energy to compress gas from near-vacuum to a given density.
These and other observations do not hold for some other pressure laws, 
which is in part why those should be considered ``exotic''.)

There are important reasons to consider flow with nonzero \defm{vorticity} $\vort=\ncurl\vv$,
which we are going to motivate by assuming \defm{irrotationality} $\vort=0$ and exploring the consequences.
Now \myeqref{eq:vsteady} can be rewritten 
\begin{eqnarray} 0 &=& \nabla\vv\dotp\vv + \nabla\piv = \nabla(\half|\vv|^2+\piv) \end{eqnarray} 
so that $\Ber$ and $\vcav$ are \emph{global} constants, same on all streamlines. 

\myeqref{eq:masssteady} yields
\begin{eqnarray} \dens\vv = -\nabla^\perp\stf \end{eqnarray} 
for a scalar function $\stf$ called \defm{stream function}. The Bernoulli relation takes the form
\begin{eqnarray} \Ber &= \half\dens^{-2}|\nabla\stf|^2 + \pif(\dens) = \rhsf(\dens,|\nabla\stf|) \end{eqnarray} 
We may use the implicit function theorem to solve for $\dens$ as long as
\begin{eqnarray} \rhsf_\dens = -\dens^{-3}|\nabla\stf|^2 + \pif_\dens(\dens) = \dens^{-3}(\ssnd^2-|\vv|^2) \end{eqnarray} 
is nonzero; it is positive as long as $|\vv|<\ssnd$ (and hence $\dens>0$), i.e.\ for \defm{subsonic} flow. 
Hence we can solve
\begin{eqnarray} \dens^{-1} = \hdiv(\half|\nabla\stf|^2) \end{eqnarray}
where $\hdiv$ is defined on some maximal interval $\cli{0}{\vmsmax}$. 

Having solved the mass and Bernoulli equations it remains to ensure irrotationality (which is needed to recover the original velocity equation \myeqref{eq:vsteady} from the Bernoulli equation):
\begin{eqnarray} 
  0 = \ncurl \vv = \ncurl \frac{-\nperp\stf}{\dens} = -\ndiv\big( \hdiv(\frac{|\nabla\stf|^2}{2}) \nabla\stf \big) 
  \label{eq:stf-divform} 
\end{eqnarray} 
After differentiation we have
\begin{eqnarray} 0 = \big(1-(\frac{v^x}{\ssnd})^2\big)\stf_{xx} - 2\frac{v^x}{\ssnd}\frac{v^y}{\ssnd}\stf_{xy} 
+ \big(1-(\frac{v^y}{\ssnd})^2\big)\stf_{yy} \label{eq:comp-potf}\end{eqnarray} 
where $\ssnd$ is a function of $\dens$, hence of $\nabla\stf$. 
The eigenvectors of the coefficient matrix $I-(\vv/\ssnd)^2$ are $\vv$ and $\vv^\perp$, 
with eigenvalues $1-\Mach^2$ and $1$ where 
\begin{eqnarray} \Mach \defeq \vlen\vv/\ssnd \end{eqnarray} 
is the \defm{Mach number}. 
Hence \eqref{eq:comp-potf} is elliptic exactly wherever it is subsonic.

The limit of decreasing Mach numbers formally (and, under various assumptions provably) yields steady irrotational 
\emph{incompressible} flow:
\begin{eqnarray} 0 &= \Lap\stf. \end{eqnarray} 
2d harmonic functions are more conveniently represented as holomorphic maps:
consider the \defm{complex velocity}
\begin{eqnarray}
  \www &\defeq& \vx-i\vy 
\end{eqnarray}
as a function of $z \defeq x+iy$. Then 
\begin{eqnarray}
  \pcz \www 
  &=&
  \half( \ndiv\vv - i\ncurl\vv ) .
\end{eqnarray}
Hence $\www$ represents an \emph{incompressible} and \emph{irrrotational} flow if and only if $\www$ is holomorphic. 
If so, it is convenient to use the \defm{complex velocity potential} $\wpot=\int^z\www~dz=\vpot+i\stf$ (which may be multivalued). 
The Cauchy-Riemann equations $\pcz\wpot=0$, together with $\pz\wpot=\www=\vx-i\vy$, yield $\vv=-\nabla\stf=\nabla\vpot$,
justifying the notation. 

At solid boundaries we use the standard \defm{slip condition}
\begin{eqnarray} 0 = \vn \dotp \vv = \vs\dotp\nabla\stf \quad, \label{eq:slip}\end{eqnarray} 
where $\vn,\vs$ are normal and tangent to the solid.
Integration along connected components of (say) a piecewise $\Cone$ boundary yields $\stf=\const$; 
if the solid boundary has a single connected component, 
then we may add an arbitrary constant to $\stf$ without changing $\vv=-\nperp\stf$ to obtain the convenient zero Dirichlet condition
\begin{eqnarray}
  \stf = 0.
\end{eqnarray}

\section{Protruding corners}

Consider a neighbourhood of a corner of a 2d solid body, as in fig.\ \myref{fig:onecorner}. 
For simplicity assume the two sides are locally straight, 
enclosing an exterior (i.e.\ fluid-side) angle $\Pola\notin\set{0,\pi,2\pi}$. 
Whatever the global problem, in many cases it is natural to seek existence of an incompressible flow by Hilbert space method, 
typically yielding an $\stf$ that is locally $\isoba1$, i.e.\ velocity $\nabla\stf$ square-integrable near the corner
with $\stf=0$ on the two sides. 
A coordinate change to $\zeta=z^{\pi/\Pola}$ maps to a small halfball; 
using $\stf=0$ on the straight side Schwarz reflection yields a harmonic function in a small ball,
with local Taylor expansion (the constant term is irrelevant for $\nabla\stf$)
\begin{eqnarray} \stf = \Im \sum_{k=1}^\infty \acf_k \zeta^k = \Im \subeq{\sum_{k=1}^\infty \acf_k z^{k\pi/\Pola}}{=\wpot(z)} \end{eqnarray}
so that
\begin{eqnarray} \www(z) = \wpot'(z) = \frac{\pi}{\Pola} \acf_1 \subeq{z^{\pi/\Pola-1}}{!} + \sum_{k=2}^\infty \frac{k\pi}{\Pola} \acf_k z^{k\pi/\Pola-1}.  \end{eqnarray} 
Here we may observe a key distinction:
if $\Pola<\pi$, i.e.\ for \emph{receding} corners, the exponent $\pi/\Pola-1$ is positive
so a square-integrable velocity $\www$ is always bounded in the corner.
This is not true at \emph{protruding} corners ($\Pola>\pi$) unless a single real scalar constraint is satisfied: $a_1=0$. 
If so, then inspection of the remaining terms shows that a nonzero $\stf$ must attain both negative and positive sign 
arbitrarily close to the corner. 
(For nonstraight but sufficiently regular sides analogous result can be obtained.)

Whereas unbounded velocity is merely undesirable for incompressible flows, 
it is mathematically \emph{impossible} in the compressible problem \myeqref{eq:comp-potf}.
The latter is quasilinear; when linearizing the operator about $\vv=0$ and $\dens=\const>0$ 
we obtain the incompressible operator $\Lap$. 
For solutions of the resulting linearized problem to be helpful in solving or understanding the nonlinear one,
we need to require unbounded velocity. 
Hence every additional protruding corner adds another scalar real constraint to the problem. 
To satisfy these constraints, a corresponding number of free real parameters may be needed!

The velocity $\www$ is holomorphic as well as bounded at infinity, so a Laurent expansion yields
\begin{eqnarray} \www(z) = c_0 + c_1z^{-1} + c_2z^{-2} + ... \myeqlab{eq:wwwww}\end{eqnarray} 
for constants $c_i$. 
$\www\conv c_0\eqdef\wwi$ at infinity (\emph{everywhere}, which is a consequence of $\vort=0$; 
rotational flows can be far more complex).
By rotational symmetry we may assume $\wwi$ is real positive; 
it is usually a chosen parameter, such as aircraft speed. 

$\Re c_1$ must be zero because it is the coefficient of a velocity term $\xx/|\xx|^2$ which would otherwise cause net mass flux
through a sufficiently large circle, 
in conflict with conservation of mass since the slip condition implies zero mass flux through the solid boundary.
Hence $c_1=\Gam/(2\pi i)$ with real 
\begin{eqnarray} 
  \Gam = \int_{\Sigma} \vv\dotp d\xx 
\end{eqnarray} 
which is called \defm{circulation}; 
$\Sigma$ is an arbitrary contour passing once counterclockwise around the circle
(due to $\vort=0$, $\Gam$ is independent of the choice of $\Sigma$, 
which may be the boundary of the body, or very large (``around infinity''), or anything in between).

Integrating $\www$ in \myeqref{eq:wwwww} yields
\begin{eqnarray}
  \wpot(z) = \wwi z + \frac{\Gam}{2\pi i} \log z + c_2 z^{-1} + ...
\end{eqnarray}
(up to an additive constant that does not change $\www=\wpot'$, hence may be taken zero). 
Consider two streamfunctions $\stf=\Im\wpot$ with same $\wwi$ and $\Gam$ and with bounded velocity. 
Then their difference $d$ is $O(z^{-1})$ at infinity, and by the slip condition converges to $0$ at the solid as well;
by the strong maximum principle it cannot have extrema in between, so $d=0$. 
Hence the flow is uniquely determined.

$\Gam$ is arbitrary for smooth bodies, as in the following explicit formula for flow around the unit circle:
\begin{eqnarray} \www(z) = \wwi ( 1 - z^{-2} ) + \frac{\Gam}{2\pi i} z^{-1}. \end{eqnarray} 
That is no longer the case if the body has protruding corners: 
if two $\stf$ with same $\wwi$ had different $\Gam$, then $d=a\log|z|+o(1)$ for $a\neq0$,
hence $\sign d=\sign a$ near infinity and thus everywhere, by the strong maximum principle and by $d=0$ at the body.
But as discussed earlier a harmonic nonzero $d$ with bounded gradient and with $d=0$ on the boundary 
must attain \emph{both} signs near the corner --- contradiction. 
Hence $\wwi$ alone uniquely determines the flow, including $\Gam$.

For a single protruding corner, and with suitable assumptions about the body, it is possible to show 
that an incompressible solution satisfying the \defm{Kutta-Joukowsky condition} of bounded velocity at the corner does exist. 
Using Blasius' theorem, Kutta-Joukowsky theory translates the resulting $\Gam$ into 
explicit and widely used formulas for the vertical force (\defm{lift}) on the body. 
At least in some physical regimes the prediction is in reasonable agreement with experimental observations 
(see fig.\ 6.7.10 and surrounding text in \cite{batchelor}). 
For particular shapes it is possible to give explicit formulas for $\www$ by conformal transformations to the unit circle
or other simple special cases.

The classical work of Frankl and Keldysh \cite{frankl-keldysh-1934}, Shiffman \cite{shiffman-exi-potf}, Bers \cite{bers-exi-uq-potf}, 
Finn-Gilbarg \cite{finn-gilbarg-uniqueness} etc.\ extends these
existence and uniqueness results and force formulas to compressible subsonic flow.

\section{Flow around bounded bodies}

The incompressible Kutta-Joukowsky theory is simple and yet powerful, 
providing insight into the all-important lift problem using only standard complex analysis. 
This combination should not obscure that irrotational inviscid flow without enhancements is inadequate for many other purposes.

\newcommand{\plateangle}{\alpha}

\begin{figure}
  \parbox{2.3in}{%
    \centerline{\input{diagonal_plate.pstex_t}}
    \caption{No subsonic solutions around non-horizontal plates}
    \label{fig:flatplate}
  }\hfil%
  \parbox{2.1in}{%
    \centerline{\input{triangle.pstex_t}}
    \caption{Utility graph argument: intersection is unavoidable}
    \label{fig:triangle}
  }%
\end{figure}

For example, what about shapes with several protruding corners? 
Since each protruding corner adds another constraint, while we still only have a single free parameter $\Gam$,
it is natural to believe that the compressible problem quickly becomes unsolvable.

Nevertheless this argument is not rigorous.
First, it is easy to find examples of nonlinear problems whose solutions cannot be obtained by simple linearization. 
Second, the constraints are not obviously independent. 
Consider for instance a flat plate (fig.\ \myref{fig:flatplate}), a degenerate example of a two-protruding-corner shape. 
If $\plateangle\in\pi\Z$ (horizontal plate), 
then compressible and incompressible problem are both obviously solvable by $\www=\const=\wwi$ in the entire domain. 

But if $\plateangle\notin\pi\Z$ (including the case of vertical plates), then we can show \cite{elling-flatplate} 
that the compressible problem does not have any uniformly subsonic solutions 
(incompressible solutions have unbounded velocity at at least one of the two corners).
We obtain such nonexistence results for a somewhat larger class of two-protruding-corner bodies;
see \cite{elling-flatplate} for details.

Due to the nonlinearity of subsonic flow, 
it is generally difficult to tell by inspection whether non-existence holds for a particular shape. 
However, using Morrey estimates \cite{elling-protrudingangle} we obtain the following rather general tool: 
if a shape does not allow incompressible solutions with bounded velocity, 
then compressible solutions with local Mach numbers below a sufficiently small positive constant do not exist either. 
In \cite{elling-polylow} we derive and apply this result for a special case of bodies with three or more protruding corners:
around an \emph{arbitrary} simple polygon there are no incompressible flows with bounded velocity, 
and therefore no low-Mach compressible flows. 

The proof of the polygon result (see \cite{elling-polylow} for full details) is based on non-planarity of the 
\defm{utility graph} (fig.\ \myref{fig:triangle}; see \cite{thomassen-1992} and references therein):
consider a bipartite graph with an edge between each of three vertices (``houses'') 
and each of the three other vertices (``utilities''); 
it is not possible to embed such a graph into the plane without allowing edges to intersect.
(Presumably the gas, water and electricity lines cannot cross because 
the utilities are uncooperative monopolies that own the land.)

Consider the sets $\set{\stf>0}$, $\set{\stf<0}$ and the interior of the body;
they are pairwise disjoint. One distinguished point in each set is chosen as ``utility''. 
Three protruding corners are chosen as ``houses''.
As discussed earlier, $\stf$ must attain both signs at each protruding corner where the velocity is locally bounded,
so all three sets must be adjacent to every house. 
In addition $\set{\stf>0}$ and $\set{\stf<0}$ are connected: both are unbounded, 
since $\stf=\wwi y+O(|z|)$ with $\wwi>0$ as $|z|\conv\infty$, 
and neither can have a bounded connected component, by the strong maximum principle. 
After choosing suitable edges through each set we would have a planar embedding of the utility graph.
The contradiction shows that any incompressible flow must have unbounded velocity at all but at most two of the protruding corners.

The proof suggests that the nonexistence is quite topological in nature. 
Simple polygons or other bodies with three or more protruding corners probably do allow irrotational flow in other settings,
for example if we restrict infinity to three ducts meeting in a junction in which the body is located. 
Such settings also have more free parameters to permit satisfying constraints at additional protruding corners. 

\begin{figure}
  \parbox{2.3in}{\input{smoothangle.pstex_t}}\hfill\parbox{2.3in}{\input{smoothanglesheet.pstex_t}}%
  \caption{Left: $\Dom$ covers an angle $\Pola$ at infinity, with variable but smooth boundaries. 
    This domain, for $\pp(\dens)=\dens^\isenc$ with $\isenc>1$,
    does not allow compressible uniformly subsonic flows whose velocity is bounded but nonvanishing at infinity. 
    Right: if vorticity is allowed, then there are trivial solutions with $\vv=\vvi=\const\neq 0$ above the vortex sheet, in particular
    near infinity at the upstream wall, whereas $\vv=0$ below the sheet.
  }
  \mylabel{fig:infangle}
\end{figure}

\section{Smooth infinite angles}

In \cite{elling-protrudingangle} we instead consider an unbounded solid:
an infinite protruding smoothened or sharp corner (see fig.\ \myref{fig:infangle} left).
It is natural to look for solutions whose velocity is bounded and converges near the upstream wall ($\graph\polahi$) 
to a prescribed nonzero constant at spatial infinity. 
Indeed for \emph{super}sonic velocity many shapes have well-known solutions based on simple waves 
\cite[section 111]{courant-friedrichs}. 
However, we prove that there do not exist any uniformly \emph{subsonic} irrotational flows of this type. 
If we do permit rotation, then some of the same shapes do allow easily constructed solutions with 
straight vortex sheet separating from the wall 
(fig.\ \myref{fig:infangle} right).

In cases where irrotational inviscid models have no solution, what models may be more appropriate?
Physical observations suggest to permit rotation, with vorticity generated at some of the protruding corners,
for example vortex sheets separated by regions of irrotational flow \cite{elling-sepsheet}. 
(Another option is to consider transonic flows, which generally feature curved shock waves that --- again --- produce vorticity.)
Vorticity is also needed for resolving part of the \defm{d'Alembert paradox} (\cite{dalembert-paradox,finn-gilbarg-uniqueness}),
namely that irrotational inviscid models predict zero drag for subsonic flows.
Although zero is a fair approximation for some shapes that are observed to have rather low drag relative to their size, 
many other shapes do have significant drag even at small viscosity. 

Historically there has been a lot of debate about how solid boundaries can generate vorticity.  
It is frequently stated that the explanation requires viscosity. 
On one hand, this is correct in the sense that physical observations --- as well as Prandtl's theory \cite{prandtl-1904} ---
confirm the essential role of thin viscous boundary layers and their instabilities. 
On the other hand, our results could be interpreted as saying that viscosity is not needed after all, 
in the sense that \emph{even in complete absence of viscosity}
some shapes and Mach numbers do not allow flows with zero vorticity.

\section*{Acknowledgements}
  This material is based upon work partially supported by the
  National Science Foundation under Grant No.\ NSF DMS-1054115.


\end{document}